\documentclass{birkjour}

\newcommand{\C}{{\mathbb C}}
\newcommand{\D}{{\mathbb D}}

\newcommand{\T}{{\mathbb T}}

\newcommand{\N}{{\mathbb N}}

\newcommand{\bmo}{{\rm BMO}}

\newcommand{\f}{\frac}
\newcommand{\ov}{\overline}

\newcommand{\ga}{\gamma}

\newcommand{\de}{\delta}
\newcommand{\la}{\lambda}
\newcommand{\ze}{\zeta}
\renewcommand{\th}{\theta}

\title[Model Subspaces in $H^1$]
{Interpolating by Functions from Model\\ 
Subspaces in $H^1$}
\author[K. M. Dyakonov]{Konstantin M. Dyakonov}

\address
{Departament de Matem\`atiques i Inform\`atica\\ 
Universitat de Barcelona, IMUB and BGSMath\\
Gran Via, 585\\ 
E-08007 Barcelona\\ 
Spain}

\address{\it \,\,\\
and}

\address{\quad\\
Instituci\'o Catalana de Recerca i Estudis Avan\c{c}ats (ICREA)\\ 
Pg. Llu\'is Companys, 23\\ 
E-08010 Barcelona\\ 
Spain}

\email{konstantin.dyakonov@icrea.cat}
\thanks{Supported in part by grants MTM2014-51834-P and MTM2017-83499-P from El Ministerio de Econom\'ia y Competitividad (Spain).}

\begin{document}
\begin{abstract}
Given an interpolating Blaschke product $B$ with zeros $\{a_j\}$, we seek to characterize the sequences of values $\{w_j\}$ for which the interpolation problem $$f(a_j)=w_j\qquad (j=1,2,\dots)$$
can be solved with a function $f$ from the model subspace $H^1\cap B\ov{H^1_0}$ of the Hardy space $H^1$.
\end{abstract}

\maketitle

\section{Introduction}

Let $\D$ stand for the disk $\{z\in\C:|z|<1\}$, $\T$ for its boundary, and $m$ for the normalized Lebesgue measure on $\T$. For $0<p\le\infty$, we denote by $H^p$ the classical Hardy $p$-space of holomorphic functions on $\D$ (see \cite[Chapter II]{G}), viewed also as a subspace of $L^p(\T,m)$. 

\par Now suppose $\th$ is an {\it inner function} on $\D$, that is, $\th\in H^\infty$ and $|\th|=1$ almost everywhere on $\T$. The associated {\it model subspace} $K^p_\th$ is then defined, this time for $1\le p\le\infty$, by putting
\begin{equation}\label{eqn:defmodsub}
K^p_\th:=H^p\cap\th\ov{H^p_0},
\end{equation}
where $H^p_0:=zH^p$ and the bar denotes complex conjugation. In other words, $K^p_\th$ is formed by those $f\in H^p$ for which the function 
$$\widetilde f:=\ov z\ov f\th$$
(living a.e. on $\T$) is in $H^p$. Moreover, the map $f\mapsto\widetilde f$ leaves $K^p_\th$ invariant. 

\par We may also rephrase the above definition by saying that $K^p_\th$ is the kernel (in $H^p$) of the Toeplitz operator with symbol $\ov\th$. For further -- and deeper -- operator theoretic connections, we refer to \cite{N}. Finally, we mention the direct sum decomposition formula 
\begin{equation}\label{eqn:dirsumdec}
H^p=K^p_\th\oplus\th H^p,\qquad1<p<\infty
\end{equation}
(with orthogonality for $p=2$), which is a consequence of the M. Riesz projection theorem; cf. \cite[Chapter III]{G}. 

\par From now on, the role of $\th$ will be played by an {\it interpolating Blaschke product}, i.e., by a function of the form $$B(z)=B_{\{a_k\}}(z):=\prod_k\f{|a_k|}{a_k}\f{a_k-z}{1-\ov a_kz},$$
with $\{a_k\}\subset\D$ and $\sum_k(1-|a_k|)<\infty$, satisfying 
$$\inf_k|B'(a_k)|\,(1-|a_k|)>0.$$
The zero sequences $\{a_k\}$ of such products are precisely the {\it interpolating sequences} for $H^\infty$ (see \cite[Chapter VII]{G}), meaning that 
\begin{equation}\label{eqn:trcarl}
H^\infty\big|_{\{a_k\}}=\ell^\infty.
\end{equation}

\par Here and throughout, we use the notation $X\big|_{\{a_k\}}$ (where $X$ is a given function space on $\D$) for the set of all sequences $\{w_k\}\subset\C$ such that the interpolation problem
\begin{equation}\label{eqn:interprob}
f(a_k)=w_k\qquad (k=1,2,\dots)
\end{equation}
has a solution $f\in X$. Equivalently, $X\big|_{\{a_k\}}$ consists of the traces $f\big|_{\{a_k\}}=\{f(a_k)\}$ that arise as $f$ ranges over $X$. 

\par Now, for $0<p<\infty$, the appropriate $H^p$ analogue of \eqref{eqn:trcarl} is known to be 
\begin{equation}\label{eqn:shashi}
H^p\big|_{\{a_k\}}=\ell^p_1(\{a_k\}),
\end{equation}
where $\ell^p_1(\{a_k\})$ stands for the set of sequences $\{w_k\}\subset\C$ with 
\begin{equation}\label{eqn:carlpsum}
\sum_k|w_k|^p\,(1-|a_k|)<\infty.
\end{equation}
Indeed, a theorem of Shapiro and Shields (see \cite{SS}) tells us that, for $1\le p<\infty$, \eqref{eqn:shashi} holds if and only if $\{a_k\}$ is an interpolating sequence for $H^\infty$; the case $0<p<1$ is, in fact, no different. Letting $B=B_{\{a_k\}}$ be the associated (interpolating) Blaschke product, we may then combine \eqref{eqn:shashi} with the observation that 
$$H^p\big|_{\{a_k\}}=K^p_B\big|_{\{a_k\}},\qquad1<p<\infty$$
(which follows upon applying \eqref{eqn:dirsumdec} with $\th=B$), to deduce that 
\begin{equation}\label{eqn:kpblp}
K^p_B\big|_{\{a_k\}}=\ell^p_1(\{a_k\}),\qquad1<p<\infty.
\end{equation}

\par The trace space $K^p_B\big|_{\{a_k\}}$ with $1<p<\infty$ is thereby nicely characterized, but the endpoints $p=1$ and $p=\infty$ present a harder problem. The latter case was recently settled by the author in \cite{DARX}. There, it was shown that a sequence $\{w_k\}\in\ell^\infty$ belongs to $K^\infty_B\big|_{\{a_k\}}$ if and only if the numbers 
\begin{equation}\label{eqn:defwtilde}
\widetilde w_k:=\sum_j\f{w_j}{B'(a_j)\cdot(1-a_j\ov a_k)}\qquad(k=1,2,\dots)
\end{equation}
satisfy $\{\widetilde w_k\}\in\ell^\infty$. We mention in passing that a similar method, coupled with \cite[Theorem 5.2]{DSpb}, yields the identity 
$$K_{*B}\big|_{\{a_k\}}=\left\{\{w_k\}\in\ell^2_1(\{a_k\}):\,
\{\widetilde w_k\}\in\ell^\infty\right\},$$
where $K_{*B}:=K^2_B\cap\bmo(\T)$. 

\par It is the other extreme (i.e., the case $p=1$) that puzzles us. We feel that the numbers $\widetilde w_k$, defined as in \eqref{eqn:defwtilde}, must again play a crucial role in the solution, and we now proceed to discuss this in detail.

\section{Problems and discussion}

\medskip\noindent\textbf{Problem 1.} Given an interpolating Blaschke product $B$ with zeros $\{a_k\}$, characterize the set $K^1_B\big|_{\{a_k\}}$. 

\medskip An equivalent formulation is as follows. 

\medskip\noindent\textbf{Problem 1$'$.} For $B=B_{\{a_k\}}$ as above, describe the functions from $K^1_B+BH^1$ in terms of their values at the $a_k$'s. 

\medskip The equivalence between the two versions is due to the obvious fact that, given a function $f\in H^1$, we have $f\in K^1_B+BH^1$ if and only if $f\big|_{\{a_k\}}=g\big|_{\{a_k\}}$ for some $g\in K^1_B$. 

\smallskip It should be noted that, whenever $B=B_{\{a_k\}}$ is an {\it infinite} interpolating Blaschke product, the trace space $K^1_B\big|_{\{a_k\}}$ is strictly smaller than $H^1\big|_{\{a_k\}}=\ell^1_1(\{a_k\})$, a result that can be deduced from \cite[Theorem 3.8]{Steg}. On the other hand, it was shown by Vinogradov in \cite{V} that 
$$K^1_{B^2}\big|_{\{a_k\}}=\ell^1_1(\{a_k\})$$ 
for each interpolating Blaschke product $B=B_{\{a_k\}}$. The subspace $K^1_{B^2}$ is, however, essentially larger than $K^1_B$. 

\par In light of our solution to the $K^\infty_B$ counterpart of Problem 1, as described in Section 1 above, the following conjecture appears to be plausible. 

\medskip\noindent\textbf{Conjecture 1.} Let $B=B_{\{a_k\}}$ be an interpolating Blaschke product. In order that a sequence of complex numbers $\{w_k\}$ belong to $K^1_B\big|_{\{a_k\}}$, it is necessary and sufficient that 
\begin{equation}\label{eqn:carlsumw}
\sum_k|w_k|\,(1-|a_k|)<\infty
\end{equation}
and 
\begin{equation}\label{eqn:carlsumwtilde}
\sum_k|\widetilde w_k|\,(1-|a_k|)<\infty
\end{equation}
(i.e., that both $\{w_k\}$ and $\{\widetilde w_k\}$ be in $\ell^1_1(\{a_k\})$).

\medskip As a matter of fact, the necessity of \eqref{eqn:carlsumw} and \eqref{eqn:carlsumwtilde} is fairly easy to verify. First of all, letting $\de_k$ denote the unit point mass at $a_k$, we know that the measure $\sum_k(1-|a_k|)\,\de_k$ is Carleson (see \cite[Chapter VII]{G}), and so 
\begin{equation}\label{eqn:carlsumfak}
\sum_k|F(a_k)|\,(1-|a_k|)<\infty
\end{equation}
for every $F\in H^1$. 
\par Now suppose that we can solve the interpolation problem \eqref{eqn:interprob} with a function $f\in K^1_B$. An application of \eqref{eqn:carlsumfak} with $F=f$ then yields \eqref{eqn:carlsumw}. Furthermore, putting $g=\widetilde f(:=\ov z\ov fB)$ a.e. on $\T$, we know that $g\in K^1_B(\subset H^1)$; in particular, $g$ extends to $\D$ by the Cauchy integral formula. Therefore, 
\begin{equation}\label{eqn:gcauchy}
\begin{aligned}
\ov{g(a_k)}&=\int_{\T}\f{\ov{g(\ze)}}{1-\ze\ov{a_k}}\,dm(\ze)\\
&=\int_{\T}\f{\ze f(\ze)\ov{B(\ze)}}{1-\ze\ov{a_k}}\,dm(\ze)\\
&=\f1{2\pi i}\int_{\T}\f{f(\ze)}{B(\ze)\cdot(1-\ze\ov{a_k})}\,d\ze,
\end{aligned}
\end{equation}
and computing the last integral by residues we find that 
\begin{equation}\label{eqn:bargak}
\ov{g(a_k)}=\sum_j\f{w_j}{B'(a_j)\cdot(1-a_j\ov a_k)}=\widetilde w_k, 
\end{equation}
for each $k\in\N$. Finally, we apply \eqref{eqn:carlsumfak} with $F=g$ and combine this with \eqref{eqn:bargak} to arrive at \eqref{eqn:carlsumwtilde}. 
\par Thus, conditions \eqref{eqn:carlsumw} and \eqref{eqn:carlsumwtilde} are indeed necessary in order that $\{w_k\}$ be in $K^1_B\big|_{\{a_k\}}$. It is the sufficiency of the two conditions that presents an open problem. 
\par In contrast to the case $1<p<\infty$, where \eqref{eqn:carlpsum} remains valid upon replacing $\{w_k\}$ by $\{\widetilde w_k\}$, no such thing is true for $p=1$. In other words, \eqref{eqn:carlsumw} no longer implies \eqref{eqn:carlsumwtilde}, so the two conditions together are actually stronger than \eqref{eqn:carlsumw} alone. Consider, as an example, the \lq\lq radial zeros" situation where the interpolating sequence $\{a_k\}$ satisfies 
\begin{equation}\label{eqn:radialcase}
0\le a_1<a_2<\dots<1
\end{equation}
(in which case the quantities $1-a_k$ must decrease exponentially, see \cite[Chapter VII]{G}). Further, let $\{\ga_k\}$ be a sequence of nonnegative numbers such that $\sum_k\ga_k<\infty$ but $\sum_k k\ga_k=\infty$, and let $w_k=B'(a_k)\cdot\ga_k$. Then 
$$(1-|a_k|)\cdot|w_k|=(1-|a_k|)\cdot|B'(a_k)|\cdot\ga_k\le\ga_k,\qquad k\in\N$$ 
(since $B$ is a unit-norm $H^\infty$ function), and \eqref{eqn:carlsumw} follows. On the other hand, 
$$\widetilde w_k=\sum_j\f{\ga_j}{1-a_ja_k}\ge\sum_{j=k}^\infty\f{\ga_j}{1-a_ja_k}
\ge\f1{2(1-a_k)}\sum_{j=k}^\infty\ga_j,$$
because for $j\ge k$ we have $a_j\ge a_k$, and hence 
$$1-a_ja_k\le1-a_k^2\le2(1-a_k).$$
Consequently, 
$$\sum_k|\widetilde w_k|\,(1-|a_k|)=\sum_k\widetilde w_k\,(1-a_k)
\ge\f12\sum_k\sum_{j=k}^\infty\ga_j=\f12\sum_j j\ga_j=\infty,$$
which means that \eqref{eqn:carlsumwtilde} breaks down. 

\par Our next problem involves the notion of an ideal sequence space, which we now recall. Suppose $\mathcal S$ is a vector space consisting of sequences of complex numbers. We say that $\mathcal S$ is {\it ideal} if, whenever $\{u_j\}\in\mathcal S$ and $\{v_j\}(\subset\C)$ is a sequence satisfying $|v_j|\le|u_j|$ for all $j$, it follows that $\{v_j\}\in\mathcal S$. Roughly speaking, this property means that the sequences belonging to $\mathcal S$ admit a nice description in terms of a certain \lq\lq size condition" on their components. 

\par When $1<p<\infty$, identity \eqref{eqn:kpblp} tells us that the trace space $K^p_B\big|_{\{a_k\}}$ coincides with $\ell^p_1(\{a_k\})$ and is therefore ideal. However, assuming that Conjecture 1 is true and arguing by analogy with the $p=\infty$ case, as studied in \cite{DARX}, we strongly believe that $K^1_B\big|_{\{a_k\}}$ is no longer ideal in general. We would even go on to conjecture that this last trace space is never ideal, as soon as $B=B_{\{a_k\}}$ is an infinite Blaschke product. The following problem arises then naturally. 

\medskip\noindent\textbf{Problem 2.} Given an interpolating Blaschke product $B$ with zeros $\{a_k\}$, determine the maximal (largest possible) ideal sequence space contained in $K^1_B\big|_{\{a_k\}}$. 

\medskip In particular, the calculations above pertaining to the special case \eqref{eqn:radialcase} show that, in this \lq\lq radial zeros" situation, the maximal ideal subspace of $K^1_B\big|_{\{a_k\}}$ must be 
$$\mathcal M:=\left\{\{w_k\}:\,\sum_kk|w_k|\,(1-a_k)<\infty\right\},$$
provided that the validity of Conjecture 1 is established. Strictly speaking, by saying that $\mathcal M$ is {\it maximal} we mean that $\mathcal M\subset K^1_B\big|_{\{a_k\}}$, while every ideal sequence space $\mathcal S$ with $\mathcal S\subset K^1_B\big|_{\{a_k\}}$ is actually contained in $\mathcal M$. 

\par Our last problem (posed in a somewhat vaguer form) deals with the map, to be denoted by $\mathcal T_{\{a_k\}}$, that takes $\{w_k\}$ to $\{\widetilde w_k\}$; here $\widetilde w_k$ is again defined by \eqref{eqn:defwtilde}, with $B=B_{\{a_k\}}$. 

\medskip\noindent\textbf{Problem 3.} Given an interpolating Blaschke product $B$ with zeros $\{a_k\}$, study the action of the linear operator 
$$\mathcal T_{\{a_k\}}:\,\{w_k\}\mapsto\{\widetilde w_k\}$$
on various (natural) sequence spaces. 

\medskip We know from \eqref{eqn:bargak} that, whenever the $w_k$'s satisfy \eqref{eqn:interprob} for some $f\in K^1_B$, the conjugates of the $\widetilde w_k$'s play a similar role for the \lq\lq partner" $\widetilde f:=\ov z\ov fB$ of $f$, so that 
$$\widetilde f(a_k)=\ov{\widetilde w_k}\qquad(k=1,2,\dots).$$
Recalling \eqref{eqn:kpblp}, we deduce that $\mathcal T_{\{a_k\}}$ maps $\ell^p_1(\{a_k\})$ with $1<p<\infty$ isomorphically onto itself. Its inverse, $\mathcal T^{-1}_{\{a_k\}}$, is then given by 
$$\mathcal T^{-1}_{\{a_k\}}:\,\{\widetilde w_k\}\mapsto
\left\{\sum_j\f{\widetilde w_j}{\ov{B'(a_j)}\cdot(1-\ov a_ja_k)}\right\}=\{w_k\},$$
which follows upon interchanging the roles of $f$ and $g(=\widetilde f)$ in \eqref{eqn:gcauchy} and \eqref{eqn:bargak}. It might be interesting to determine further sequence spaces that are preserved by $\mathcal T_{\{a_k\}}$. Anyhow, this $\mathcal T_{\{a_k\}}$-invariance property can no longer be guaranteed for the endpoint spaces $\ell^1_1(\{a_k\})$ and $\ell^\infty$. In fact, we have seen in \cite{DARX} that the set of those $\ell^\infty$ sequences that are mapped by $\mathcal T_{\{a_k\}}$ into $\ell^\infty$ coincides with the trace space $K^\infty_B\big|_{\{a_k\}}$ (which is, typically, smaller than $\ell^\infty$); and Conjecture 1 says that at the other extreme, when $p=1$, the situation should be similar. 
\par Furthermore, the suggested analysis of the $\mathcal T_{\{a_k\}}$ operator could perhaps be extended, at least in part, to the case of a generic (non-interpolating) Blaschke sequence $\{a_k\}$. 
\par Finally, going back to Problem 1, we remark that the range $0<p<1$ is also worth looking at, in addition to the $p=1$ case. It seems reasonable, though, to define the corresponding $K^p_B$ space as the closure of $K^\infty_B$ in $H^p$, rather than use the original definition \eqref{eqn:defmodsub} for $p<1$. The problem would then consist in describing the trace space $K^p_B\big|_{\{a_k\}}$ for an interpolating Blaschke product $B=B_{\{a_k\}}$. More generally, one might wish to characterize $K^p_\th\big|_{\{\la_k\}}$, say with $0<p\le1$ and $\th$ inner, for an arbitrary  interpolating sequence $\{\la_k\}$ in $\D$. This time, however, the chances for a neat solution appear to be slim. Indeed, even the \lq\lq good" case $1<p<\infty$ is far from being completely understood; see, e.g., \cite{DPAMS, HNP} for a discussion of such interpolation problems in $K^2_\th$.

\end{document}